\documentclass[letterpaper, 10 pt, conference]{ieeeconf}

\IEEEoverridecommandlockouts                              

\usepackage[linesnumbered,ruled]{algorithm2e}
\usepackage{amsmath} 
\usepackage{amssymb}  
\usepackage[normalem]{ulem} 
\usepackage{bm}
\usepackage{xcolor}
\usepackage{comment}
\usepackage[hidelinks]{hyperref}
\usepackage[font=scriptsize]{caption}
\usepackage{graphicx} 
\usepackage{svg}
\usepackage{import}
\usepackage[font=scriptsize]{subfig}
\usepackage{siunitx}
\usepackage{soul}

\newtheorem{theorem}{Theorem}

\newtheorem{assumption}{Assumption}

\newtheorem{remark}{Remark}

\usepackage[maxbibnames=10,giveninits=true,sorting=none]{biblatex} 
\AtEveryBibitem{%
  \clearfield{issn} 
  \clearfield{doi} 
  \clearfield{url}
  \clearfield{editor}
}
\DeclareFieldFormat*{year}{{}}
\renewbibmacro{in:}{} 

\newcommand{\tr}[1]{{#1}^{\ensuremath{\mathsf{T}}}} 

\usepackage{xcolor}

\title{Convexity in Optimal Control Problems}

\author{Abhijeet, Mohamed Naveed Gul Mohamed, Aayushman Sharma, and Suman Chakravorty  
\thanks{The authors are with the Department of Aerospace Engineering, Texas A\&M University, College Station, TX 77843 USA. \{\tt abhinir, naveed, aayushmansharma, schakrav\}@tamu.edu
}}

\begin{document}

\maketitle

\begin{abstract}
This paper investigates the central role played by the Hamiltonian in continuous-time nonlinear optimal control problems. We show that the strict convexity of the Hamiltonian in the control variable is a sufficient condition for the existence of a unique optimal trajectory, and the nonlinearity/non-convexity of the dynamics and the cost are immaterial. The analysis is extended to discrete-time problems, revealing that discretization destroys the convex Hamiltonian structure, leading to multiple spurious optima, unless the time discretization is sufficiently small. We present simulated results comparing the ``indirect" Iterative Linear Quadratic Regulator (iLQR) and the ``direct" Sequential Quadratic Programming (SQP) approach for solving the optimal control problem for the cartpole and pendulum models to validate the theoretical analysis. Results show that the ILQR always converges to the ``globally" optimum solution while the SQP approach gets stuck in spurious minima given multiple random initial guesses for a time discretization that is insufficiently small, while both converge to the same unique solution if the discretization is sufficiently small.
\end{abstract}

\section{\uppercase{Introduction}}
Optimal nonlinear control has been the subject of extensive research and application for well over half a century \cite{lewis2012optimal, bryson}. The solution of a nonlinear optimal control problem can be obtained by solving the necessary conditions obtained from Pontryagin's minimum principle and involves the challenging task of solving a so-called two-point boundary value problem (TPBVP) \cite{bryson,kirk2004optimal} which leads to an open-loop solution. Alternatively, the optimal control problem can be posed as a Dynamic Programming (DP) problem in discrete time; equivalently, the Hamilton-Jacobi-Bellman (HJB) equation in continuous time, which, in theory, when solved, provides a feedback solution to the optimal control problem \cite{bellman1966dynamic,bertsekas2012dynamic}. However, the DP problem is known to be subject to the ``Curse of Dimensionality," making it computationally intractable for most systems. Thus, the open-loop approach of solving a TPBVP is computationally attractive, and one may recover feedback by a perturbation feedback \cite{bryson,naveedACC,mohamed2023feedback} or recursively replanning as in Model Predictive Control (MPC) \cite{Mayne_1}. 
The TPBVPs specify initial conditions on the state and terminal conditions on the so-called co-state. However, the nonlinearity of the equations complicates the solution, potentially leading to multiple solutions. 
This leads us to the primary focus of the paper, where we demonstrate that if the Hamiltonian of a continuous-time optimal control problem is strictly convex in the control variable, then a unique global solution to the continuous-time optimal control exists and is given by the solution to the two-point boundary value problem, i.e., the necessary conditions for optimality are also sufficient in such a case. This implies that the nonlinearity/ non-convexity of the cost and dynamics are immaterial to the existence of a unique global minimum to the optimal control problem; rather, it is their interplay via the Hamiltonian that is critical. \\

However, the presence of a unique solution does not necessarily imply that the problem can be easily solved. This leads us to the second major area of focus of the paper. Having established the existence of a unique solution, our inquiry delves into the alternate methods for its determination. Most optimal control problem solutions involve the time discretization of the problem and employ suitable numerical techniques to obtain a solution for the resulting discretization problem. The so-called ``direct"  methods have emerged as the most widely used technique over the past two decades for solving these problems, offering a general-purpose framework for handling constraints during problem-solving. The underlying idea behind these solvers is to discretize the optimal control problem using computationally efficient pseudo-spectral methods and solve the resulting discrete Nonlinear Programming (NLP) problem \cite{RF1,GPOPS2,casadi}. It is common to employ general-purpose NLP solvers to solve the discretized problems  \cite{snopt,ipopt}.\\
An alternative is the so-called ``indirect approach," wherein one tries to solve the necessary conditions for the optimal control problem without recourse to an NLP solver. Typically, this involves an Euler discretization of the optimal control problem, which is then solved using a suitable computational approach. Notable among such approaches is the differential dynamic programming (DDP) \cite{DDP} and its closely allied cousin, the Iterative Linear Quadratic Regulator (iLQR) \cite{ilqg,ILQG_tassa2012synthesis}. The main difference between the approaches is that DDP expands the nonlinear dynamics to the second order while the ILQR only expands to the first order. Nonetheless, the performance of the approaches is identical. The drawback of the DDP/ ILQR approach is that they work for problems without state/ control constraints as opposed to NLP-based approaches; however, their great advantage is that they exploit the structure of the optimal control problem and are far more efficient when compared to NLP approaches, their time complexity being $O(N)$ versus $O(N^3)$ for NLP, where $N$ is the time horizon of the problem.\\
It should be noted that, in order to solve it, we convert a continuous-time problem to a discrete-time problem. However, it is important to recognize that solving a discrete-time problem differs from solving its continuous-time counterpart. Thus, we explore and elucidate how converting a problem to discrete-time can substantially alter the problem, potentially leading to spurious minima. We show that the continuous time convex Hamiltonian optimal control structure breaks down if the discretization time is not sufficiently small. In particular, we demonstrate that discretization with an insufficiently small time step may yield multiple solutions to a continuous time problem that was initially expected to possess a unique globally optimum solution due to the convex Hamiltonian structure. 

Furthermore, we compare the performance of NLP solvers with that of the iLQR. We conduct computational experiments using both methods and observe that iLQR tends to converge to the same solution regardless of the initial guess. On the other hand, NLP solvers return multiple different solutions to the discrete-time problem when initiated with different initial guesses. We also demonstrate that both methods converge to a single unique solution if the time discretization is made small enough.\\

The rest of the paper is organized as follows. In Section \ref{sec:2}, we introduce the problem and the assumptions underlying the analysis. In Section \ref{sec:3}, we present our central results regarding the structure of the Hamiltonian in continuous and discrete time. Section \ref{sec:comp_meth} is dedicated to outlining the ``direct" SQP  and the ``indirect" iLQR solvers for the optimal control problem. Finally, we present computational results for the optimal control of the pendulum and cartpole systems in Section \ref{sec:empirical_Results}.

\section{\uppercase{Preliminaries}}
\label{sec:2}
In the following, we detail the optimal control problem that we seek to solve and the assumptions underlying the problem.
\subsection{Optimal Control Problem}
In this section, we outline the concepts and formulation which will be used in the next sections. We are solving a nonlinear optimal control problem here, and we comment on the nature of the solution obtained for a particular class of problems. We consider the minimization of an objective function of the form
\begin{align}\label{eq:cost_function_cont}
    \min_{\mathbf{u}} \mathcal{J} = \Phi_{t_{f}}(\mathbf{x}_{t_{f}}) + \int_{0}^{t_{f}} [\mathcal{L}(\mathbf{x}) + \mathcal{R}(\mathbf{u})] dt,
\end{align}
where $\Phi_{t_{f}}(\mathbf{x}_{t_{f}})$ is the terminal cost, and $\mathcal{L}(\mathbf{x}) + \mathcal{R}(\mathbf{u})$ is the Lagrangian/ incremental cost of the problem. The optimal control problem is subject to control-affine dynamic constraints of the form 
\begin{align}\label{eq:state_eq_cont}
    \dot{\mathbf{x}} = \mathcal{F}(\mathbf{x}) + \mathcal{G}(\mathbf{x})\mathbf{u},
\end{align}
where the state vector $\mathbf{x}$ $\in$ $\mathbb{R}^{n}$, the control input $\mathbf{u}$ $\in$ $\mathbb{R}^{m}$, $\mathcal{F}(\mathbf{x})$ = $\begin{bmatrix}
    f_{1}(\mathbf{x})\\
    f_{2}(\mathbf{x})\\
    \vdots\\
    f_{n}(\mathbf{x})
\end{bmatrix}$, and $\mathcal{G}(\mathbf{x})$ is $\begin{bmatrix}
    \Gamma_{11}(\mathbf{x})&& \cdots && \Gamma_{1m}(\mathbf{x})\\
    \Gamma_{21}(\mathbf{x}) && \cdots && \Gamma_{2m}(\mathbf{x})\\
    \vdots && \ddots && \vdots \\
    \Gamma_{n1}(\mathbf{x})&& \cdots && \Gamma_{nm}(\mathbf{x}) 
\end{bmatrix}$. $\Gamma_{ij}(\mathbf{x})$ are nonlinear functions of the state $\mathbf{x}$.

\subsection{Assumptions}
We make the following assumptions on the functions involved in the optimal control problem.
\begin{assumption}{(A1)} \label{assump:1}
\textit{Smoothness.} We assume that the functions $f_{i}(\mathbf{x})$, $\Gamma_{ij}(\mathbf{x})$, $\mathcal{L}(\mathbf{x})$, and $\Phi_{t_{f}}(\mathbf{x}_{t_{f}})$ are all twice continuously differentiable ($C^{2}$) functions.
\end{assumption}

\begin{assumption}{(A2)} \label{assump:2}
\textit{Convexity of the cost.} We assume that the control cost function $\mathcal{R}(\mathbf{u})$  is a strict convex function with respect to the control inputs, $\mathbf{u}$.
\end{assumption}



\section{\uppercase{Convexity in Optimal Control}}\label{sec:3}
In this section, we show that there exists a unique solution for the optimal control problem in continuous time under the assumptions made in the previous section which render the Hamiltonian strictly convex in the control variable. Next, we show that this convex Hamiltonian structure breaks down when the problem is discretized unless the discretization time is small enough, underscoring the importance of a fine enough discretization.

\subsection{Existence of a Unique Solution for the Continuous-time Nonlinear Optimal Control Problem}
We begin with a nonlinear optimal control problem and show that the convexity of the Hamiltonian is a sufficient condition for the existence of a unique solution.\\ 
Consider the optimal control problem described by Eqs. \eqref{eq:cost_function_cont} and \eqref{eq:state_eq_cont}. For comprehensibility, we will develop and show our results for the scalar case. The generalization to the vector case is straightforward, albeit somewhat notationally tedious. Consider the objective function
\begin{equation}
    \min_{u} J = \phi(x_{t_{f}}) + \int_{0}^{t_{f}} [l(x) +  R(u)]dt,
    \label{eq:cost_function_scalar}
\end{equation}
subject to one-dimensional dynamic constraints 
\begin{equation}
    \dot{x} = f(x) + g(x)u.
    \label{eq:state_eq_scalar}
\end{equation}

For the optimal cost represented in Eq.(\ref{eq:cost_function_scalar}), the Hamiltonian can be written as
\begin{align}
    \mathcal{H} = l(x) + R(u) + \lambda[f(x) + g(x)u].  
\end{align}

It follows from the first order necessary condition that
\begin{align}
    \frac{\partial \mathcal{H}}{\partial u} = 0, \hspace{3mm}\dot{\lambda} = -\frac{\partial \mathcal{H}}{\partial x}, and \hspace{3mm} \dot{x} =\frac{\partial \mathcal{H}}{\partial \lambda}.
\end{align}
\begin{align}
    \frac{\partial \mathcal{H}}{\partial u} = 0 \Rightarrow \frac{\partial R(u)}{\partial u} + \lambda g(x) = 0
\end{align}

\begin{align}
    \Rightarrow R_{u}(u) + \lambda g(x) = 0,
    \label{eq:con_result}
\end{align}
where $R_{u} = \frac{\partial R}{\partial u}$. Here, $R_{u}^{-1}$ (inverse of the function $R_u(.)$) always exists and the corresponding solution $u^{*}$ is unique since $R$ is strictly convex \cite{boyd2004convex}. Thus,
\begin{align}
    u^{*} = R_{u}^{-1}(\lambda g(x)).
\equiv \Psi(\lambda g(x)),
    \label{eq:control_val}
\end{align}
where we have replaced the function $R_{u}^{-1}$ with $\Psi$. In addition to this, we have
\begin{align}
    \dot{\lambda} = -l_{x} - \lambda[f_{x} + g_{x}\Psi(\lambda g)]\\
    \dot{x} = f + g\Psi(\lambda g)
    \label{eq:state_tpbvp}
\end{align}
where we have dropped the arguments of $f$, $g$ and $l$ for sake of notation simplicity. With the initial condition $x(0)$ = $x_{0}$ and terminal condition $\lambda(t_{f})$ = $\phi_{x}(x_{t_{f}})$, this becomes a two-point boundary value problem (TPBVP). 
The optimal control problem is equivalent to solving the Hamilton-Jacobi-Bellman (HJB) partial differential equation, which can be written as:
\begin{align}
    -\frac{\partial J}{\partial t} = \min_{u}[\mathcal{H}],
\end{align}
where the co-state having the interpretation $\lambda = \frac{\partial J}{\partial x}.$
The HJB equation is solved backwards in time from the terminal condition $J(t_{f},x) = \phi(x_{t_{f}})$. Define $p = \frac{\partial J}{\partial t}$, then we can write the HJB as:
\begin{align}
    F(t,x,J,p,\lambda) = p + l + R(\Psi(\lambda g)) + \lambda[f + g\Psi(\lambda g)] = 0,
\end{align}
The HJB is a first-order nonlinear PDE and, consequently, can be solved using the Method of Characteristics \cite{Courant-Hilbert}.
Using the Method of Characteristics, we can turn the PDE into a family of ODEs \cite{naveedACC}, and these ODEs can be written as: $\dot{x}$ = $F_{\lambda}$ and $\dot{\lambda}$ = $-F_{x} - \lambda F_{J}$.

\begin{align*}
    \dot{x} = F_{\lambda} = \Bigg{[}\frac{\partial R(\Psi(\lambda g))}{\partial \Psi(\lambda g)}\Bigg{]} \Bigg{[}\frac{\partial \Psi(\lambda g)}{\partial \lambda}\Bigg{]} +\\ f + g\Psi(\lambda g) + \lambda g\Bigg{[}\frac{\partial \Psi(\lambda g)}{\partial \lambda }\Bigg{]},
\end{align*}
Using Eq. \eqref{eq:con_result},
\begin{align}
    \dot{x} &= -\lambda g \Bigg{[}\frac{\partial \Psi(\lambda g)}{\partial \lambda}\Bigg{]} + f + g\Psi(\lambda g) + \lambda g\Bigg{[}\frac{\partial \Psi(\lambda g)}{\partial \lambda }\Bigg{]} \nonumber \\
    \Rightarrow \dot{x} &= f + g\Psi(\lambda g).
    \label{eq:state_eq_con_subs}
\end{align}
The result obtained in Eq. \eqref{eq:state_eq_con_subs} is sensible as it gives us the closed-loop state equations, which could have been obtained by substituting Eq. \eqref{eq:control_val} in Eq. \eqref{eq:state_eq_scalar}. We further have the co-state equations as
\begin{align*}
    \dot{\lambda} = -F_{x} - \lambda \underbrace{F_{J}}_{= 0},
\end{align*}

\begin{align*}
    \dot{\lambda} = -l_{x} - \Bigg{[}\frac{\partial R(\Psi(\lambda g))}{\partial \Psi(\lambda g)}\Bigg{]} \Bigg{[}\frac{\partial \Psi(\lambda g)}{\partial (\lambda g)}\Bigg{]}\lambda g_{x}   -\\
    \lambda f_{x} - \lambda g_{x}\Psi(\lambda g) - \lambda g \Bigg{[}\frac{\partial \Psi(\lambda g)}{\partial (\lambda g)}\Bigg{]}\lambda g_{x},
\end{align*}

\begin{align*}
    \dot{\lambda} = -l_{x} + \lambda g \Bigg{[}\frac{\partial \Psi(\lambda g)}{\partial (\lambda g)}\Bigg{]}\lambda g_{x}   -
    \lambda f_{x} -\\ \lambda g_{x}\Psi(\lambda g) - \lambda g \Bigg{[}\frac{\partial \Psi(\lambda g)}{\partial (\lambda g)}\Bigg{]}\lambda g_{x},
\end{align*}

\begin{align}
    \dot{\lambda} = -l_{x} - \lambda f_{x} - \lambda g_{x} \Psi(\lambda g).
    \label{eq:costate_eq_con_subs}
\end{align}
The equations presented above need to be solved backwards in time, starting from the terminal condition $(x_{t_f}, \lambda(t_f) = \phi_x(x_{t_f}))$, to derive the characteristic curves of the Hamilton-Jacobi-Bellman (HJB) equation for various terminal states $x_{t_f}$. Note that the state and co-state equations are the same as in the TPBVP, except that the characteristic equations are solved backwards in time.\\
Assuming the optimal trajectory for the TPBVP is defined by $u^{*}(t)$ and $x^{*}(t)$, if we use the characteristic equations to backward propagate from the terminal conditions $x^{}(t_{f})$ and $\lambda^{*}(t_{f}) = \phi_{x}(x^{*}(t_{f}))$, the resulting curve passes through $x_0$. Next, we leverage a crucial finding established in \cite{mohamed2023feedback}, which states that at any time $t$, the co-state $\lambda(t)$ can be uniquely expressed as a function $\psi^{t}(x_{t})$ of the state $x_t$ owing to the uniqueness of the solutions of the Lagrange-Charpit equations, assuming that the equations have a solution on the interval $[0,t_f]$ \cite[Proposition~4]{mohamed2023feedback}.
Since the characteristic curve backward in time passes through $x_{0}$ starting at  the terminal state $x_{t_{f}}^{*}$, along with $\lambda^{*}(t_{f})$  and satisfies the Lagrange-Charpit equations (Eq. \eqref{eq:state_eq_con_subs} and \eqref{eq:costate_eq_con_subs}), it follows that the initial co-state, $\lambda_0$, is uniquely determined by $x_{0}$. Moreover, this solution adheres to the minimum principle as it constitutes a solution to the original two-point boundary value problem by definition. It's apparent that since the solution to the necessary conditions is unique, it also represents the ``global optimum". Thus, in the development above, we have proved the following result.

\begin{theorem}
Let the Lagrange-Charpit equations (Eq. \eqref{eq:state_eq_con_subs} and \eqref{eq:costate_eq_con_subs}) have a solution on $[0,t_f]$ for any $x_0$. Then, for any initial state $x_0$, under A\ref{assump:1} and A\ref{assump:2}, there exists a unique global minimum solution to the optimal control problem specified by Eq. \eqref{eq:cost_function_cont} and Eq. \eqref{eq:state_eq_cont} which can be found by solving the TPBVP specified by Eq. \eqref{eq:control_val}-\eqref{eq:state_tpbvp}.
\end{theorem}
\begin{remark}
    Thus, we can assert that the strict convexity of the Hamiltonian in the control variable is a sufficient condition for ensuring the existence of a unique solution to the TPBVP underlying the necessary conditions. Hence, the nonlinearity/ nonconvexity of the cost and dynamics are immaterial to the uniqueness of the solution; rather, it is the interaction of the two via the Hamiltonian that is critical. In fact, the underlying reason is that an optimal control problem is not simply an optimization problem but also that the necessary conditions constitute a closed-loop dynamical system. It is precisely the HJB-based dynamical system interpretation that allows us to conclude that the solution is unique.
\end{remark}

\subsection{The Discrete-Time Case}\label{subsec:discrete-time}
In most applications, we have to discretize the system. Notably, virtually all techniques employed in the field of optimal control rely on some form of discretization to tackle continuous-time problems. 
Again, limiting to scalar systems for simplicity, the discrete-time system dynamics are given by:
\begin{align}
    x_{k} = q(x_{k-1}, u_{k-1}),
\end{align}
and the discrete-time objective becomes (with a forward Euler approximation):
\begin{align}
    \min_{x_{k},u_{k}} J = \Phi(x_{N}) + \sum_{i=0}^{N-1}[\mathcal{L}(x_{k}) + \mathcal{R}(u_{k})]\Delta t,
    \label{eq:discrete_cost}
\end{align}
where $N=T/\Delta t$. However, note that the state:
\begin{align}
    x_k = x_{k-1} + \int_{(k-1)\Delta t}^{k \Delta t} f(x(\tau)) + g(x(\tau)) u_{k-1}
    \label{eq:disrete_dynamics}
d\tau = \nonumber\\
x_{k-1} + \int_{(k-1)\Delta t}^{k\Delta t} f(x(\tau)) d\tau + \Bigg{(}\int_{(k-1)\Delta t}^{k\Delta t} g(x(\tau)) d\tau\Bigg{)} u_{k-1},
\end{align}
and, in general, $x(\tau)$ is a nonlinear function of both the state $x_{k-1}$ and the control $u_{k-1}$, and hence, the discrete time state dynamics $q(x_{k-1}, u_{k-1})$ is not control affine. The affine control structure remains only if the sampling time is sufficiently small so that a forward Euler approximation holds in between the time steps.
Next, we write the Dynamic Programming equation for the evolution of the optimal cost function for the discrete-time system above ( the HJB is the continuous-time version of the DP equation). Let $c(x,u) = (\mathcal{L} (x) + \mathcal{R}(u)) \Delta t$, then the DP equation is:
\begin{align}
    J_k(x) = \min_u [Q_{k+1}(x,u)], \nonumber\\
    Q_{k+1}(x,u) = c(x,u) + J_{k+1}(q(x,u)),
\end{align}
where note that the $Q$-function plays the role of the Hamiltonian $\mathcal{H}$ in continuous time. The solution above is swept back from the terminal cost function $J_N(x) = \Phi(x)$. However, note now that owing to the loss of the control affine structure in the discrete time dynamics, the Q-function does not have a unique minimizer, in general. In fact, there may be multiple minima of the Q-function. Thus, in general, it is not possible to write down a difference equation analogous to the characteristic differential equations that were written for the HJB in continuous time. Hence, in the discrete-time case, the lack of strict convexity of the Q-function implies that there may be multiple minima to the problem, as we shall show in our empirical results shortly. Thus, maintenance of the control affine structure int he dynamics is critical in obtaining a global optimum solution, which can only be assured in discrete time if the sampling time is taken to be sufficiently small so as to be a faithful approximation of the continuous time problem.

\section{Computational Approach to Optimal Control}\label{sec:comp_meth}
In this section, we outline two of the most widely used approaches to solve the nonlinear optimal control problem presented in the paper. We discuss the ``direct" approach, where the discretized optimal control problem is treated as an NLP and solved using a general-purpose NLP solver such as Sequential Quadratic Programming (SQP), alongside the ``indirect" approach, where the necessary conditions of optimality for the discretized optimal control problem are satisfied utilizing the inherent causal structure of control resulting in the Iterative Linear Quadratic Regulator (iLQR). 

\subsection{Direct Approach: Discrete-time Sequential Quadratic Programming}

Sequential Quadratic programming stands out as one of the most widely used nonlinear programming techniques to solve discrete-time optimal control problems. We present the core structure of the algorithm to juxtapose it with iLQR. For the discrete-time problem mentioned in subsection \ref{subsec:discrete-time}, the objective is given by Eq. \eqref{eq:discrete_cost} subject to \eqref{eq:disrete_dynamics}. SQP solves this problem iteratively by forming a linear approximation of the constraints and a quadratic approximation of the cost around the current estimate of the solution, and solving the resulting quadratic program (QP) to find a descent direction for the optimization. This transition introduces the perturbations $\delta x_{k}$ and $\delta u_{k}$  the current estimate and results in the QP:
\begin{align}
    &\min_{\delta x_{k},\delta u_{k}} \delta J = \Phi_{x}(x_{N}) \delta x_{N} + \frac{1}{2}\tr{\delta x_{N}}\Phi_{xx}(x_{N})\delta x_{N} +\nonumber\\
    &\hspace{20mm} \sum_{i=0}^{N-1}[\mathcal{L}_{x}\delta x_{k} + \frac{1}{2}\tr{\delta x_{k}}\mathcal{L}_{xx}\delta x_{k}  +\nonumber\\ &\hspace{20mm} \mathcal{R}_{u}\delta u_{k} + \frac{1}{2}\tr{\delta u_{k}}\mathcal{R}_{uu}\delta u_{k}]\Delta t,\nonumber \\
    &s.t. \hspace{4mm} \zeta_{k}(\delta x_{k}, \delta x_{k-1},\delta u_{k-1}) =\nonumber \\ & \hspace{20mm}\delta x_{k} - A_{k-1}\delta x_{k-1} - B_{k - 1}\delta u_{k-1} = 0,
    \label{eq:sqp_cost}
\end{align}
where all the parameters of the QP are found from expanding the cost to the second order and the dynamics to the first order, around the current solution.
It should be noted here that we have neglected some higher order terms and omitted the terms not involving $\delta x_{k}$ or $\delta u_{k}$ in Eq. \eqref{eq:sqp_cost} as they are constants. We initialize values of $x_{k}$ and $u_{k}$. Then, we recursively update it using values of $(\delta x_{k},\delta u_{k})$ until convergence. 

\subsection{Indirect Approach: Iterative Linear Quadratic Regulator}
The iLQR is identical to the SQP as it iteratively linearizes the dynamical constraints, quadratizes the cost around the current estimate and recursively updates the solution to reach the optimum. However, the ILQR recognizes that every successive QP is an LQR problem, and thus, the solution to the component QPs is found by a highly efficient forward-backwards sweep that utilizes the causal structure of the optimal control problem, as shown below.
The iLQR algorithm entails guessing an initial $u_{k}$ and updating the state trajectory $x_{k}$ using this guessed value. Subsequently, the value of $u_{k}$ is updated as $u_{k}^{\eta +1} = u_{k}^{\eta} + \delta u_{k}$, where $\eta$ denotes the iteration number, and $\delta u_{k}$ is given by
\begin{align}
    \delta u_{k} = -k_{k} - K_{k}\delta x_{k},
\end{align}
where $k_{k}$ and $K_{k}$ are parameters updated during backward propagation. In iLQR, the co-state $\lambda_{k}$ has the feedback form $- v_{k} - V_{k}\delta x_{k}$ and the first-order necessary condition gives us a recursive update on $v$ and $V$ given by
\begin{align}
    v_{k} = -\mathcal{L}_{x} - \tr{A_{k}} v_{k+1} +  \tr{A_{k}}V_{k+1}B_{k} \nonumber \\ (\mathcal{R}_{uu} + \tr{B_{k}}V_{k+1}B_{k})^{-1}
    (\tr{B}v_{k+1} + \mathcal{R}_{uu}u_{k}^{n}) \nonumber,\\
    V_{k} = \mathcal{L}_{xx} + \tr{A_{k}}(V_{k+1}^{-1} + B_{k}(\mathcal{R}_{uu}(u_{k}^{n}))^{-1}\tr{B_{k}})^{-1}A_{k},
\end{align}
where we have terminal constraints $V_{N}$ = $\phi_{xx}$, and $v_{N}$ = $\phi_{x}x_{N}$. These can be used to obtain
\begin{align}
    k_{k} = (\mathcal{R}_{uu} + \tr{B_{k}}V_{k+1}B_{k})^{-1}(\mathcal{R}_{uu}u_{k} + \tr{B_{k}}v_{k+1}), \nonumber\\
    K_{k} = (\mathcal{R}_{uu} + \tr{B}V_{k+1}B)^{-1}\tr{B_{k}}V_{k+1}A_{k}.
\end{align}
This sums up the steps involved in calculating the optimal trajectory using iLQR.\\
Under assumptions A\ref{assump:1}, A\ref{assump:2} and a sufficiently small discretization, it can be shown that the assumed optimal control problem is guaranteed to converge to a global minimum when using iLQR \cite{wang2022search1, ranwangCDC}.

\begin{remark}
    The number of variables in the iterative QP problem is $(n+m)N$, where we note that we solve for both the state and control for every time step. In an NLP solver such as SQP, all the structure of the problem is lost, and thus, the complexity of the computation is $O((n+m)^3N^3)$ as solving the necessary conditions at every recursion amounts to solving a linear equation in $(n+m)N$ variables. However, since ILQR utilizes the causal structure inherent in the control problem, it has complexity $O((n+m)^3 N)$, which is a tremendous saving when the time horizons involved are long.
\end{remark}
\begin{remark}{\label{remark3}}
    As we shall see in the following Section, ILQR typically converges to a single solution with different random initial guesses for the control input, regardless of the selection of $\Delta t$. In contrast, SQP converges to various spurious minima. The plausible reason may be the structure inherent in the iLQR algorithm, optimization over a lesser number of variables, and automatic satisfaction of the dynamic constraints.
\end{remark}

\section{\uppercase{Empirical Results}}\label{sec:empirical_Results}
In this section, we provide empirical evidence to corroborate our theory. We consider a continuous-time nonlinear optimal control problem and discretize the problem in order to it. We demonstrate that the discrete-time problem converges to various different solutions, starting from different initial guesses, if the time discretization is not small. We also show that the solution does converge to a unique solution as we reduce the time step to a fine size. We compare the NLP and ILQR solutions against each other, keeping in mind that for a fine enough discretization, there should be a unique solution. \\
Due to the paucity of space, we consider two dynamical systems: the pendulum and the cart pole. For the pendulum, the length of the rod is taken to be 0.5m, and the mass is 0.5 Kg. For the cart pole, the mass of the cart is 1 Kg, and the length and the mass of the pole are 0.6m and 0.01 Kg, respectively.  
\subsection{Continuous-time Optimal Control Problem with Convex Hamiltonian}
We consider quadratic cost functions for both the pendulum and the cartpole. The quadratic cost makes the Hamiltonian convex. The quadratic cost is given by
\begin{equation}
    J = \frac{1}{2}\tr{x_{t_{f}}}S_{f}x_{t_{f}} + \int_{0}^{t_{f}} [\tr{x}Qx + \tr{u}Ru] dt,
\end{equation}

where $S_{f}$, $Q$, and $R$ are chosen to be diagonal matrices with some weight. The subscript $t_{f}$ denotes the value of state calculated at the final/terminal time. The values of $Q$, $R$ and $S_{f}$ for both cartpole and pendulum are chosen to be $100I_{4 \times 4}$, $10$, and $1000I_{4 \times 4}$, where $I$ represents the identity matrix.

\subsection{Sequential Quadratic programming (SQP): Convergence to Spurious Minima}
\begin{figure}[!htbp]
    \centering
   \subfloat[Solution 1 ($J = 6.103895e+03$).]{\includegraphics[width=.49\linewidth]{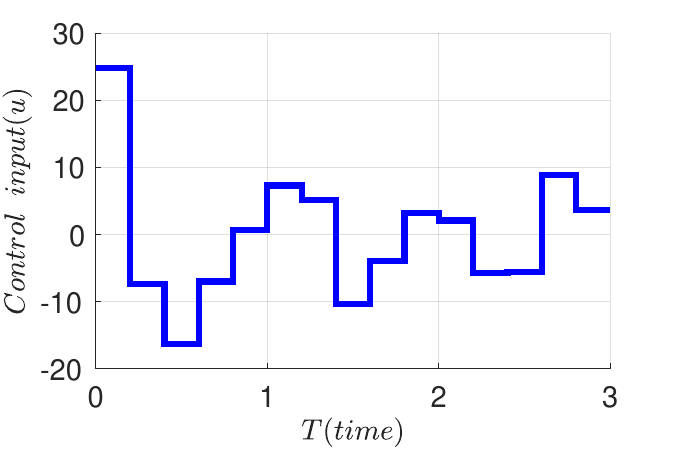}}
   \hfill
    \subfloat[Solution 2 ($J = 1.650793e+04$).] {\includegraphics[width=0.49\linewidth]{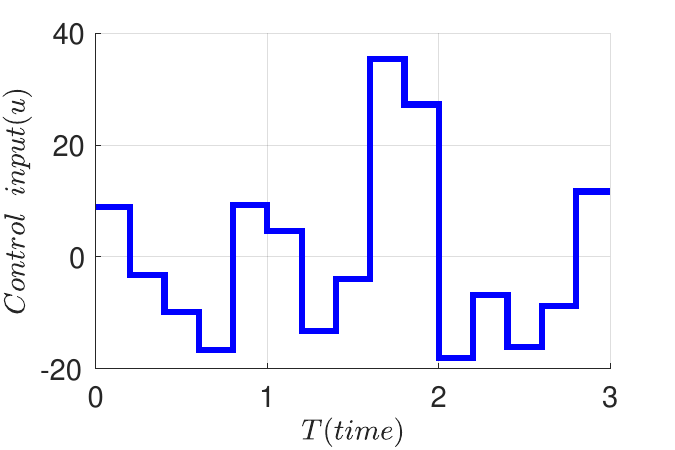}}%
    \hfill
    \subfloat[Solution 3 ($J = 1.941844e+03$).] {\includegraphics[width=0.49\linewidth]{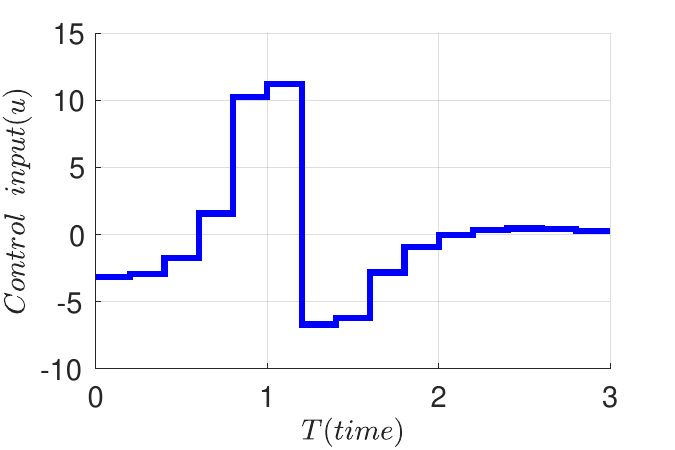}}
    \hfill
    \subfloat[Solution 4 ($J = 2.052650e+04$).] {\includegraphics[width=0.49\linewidth]{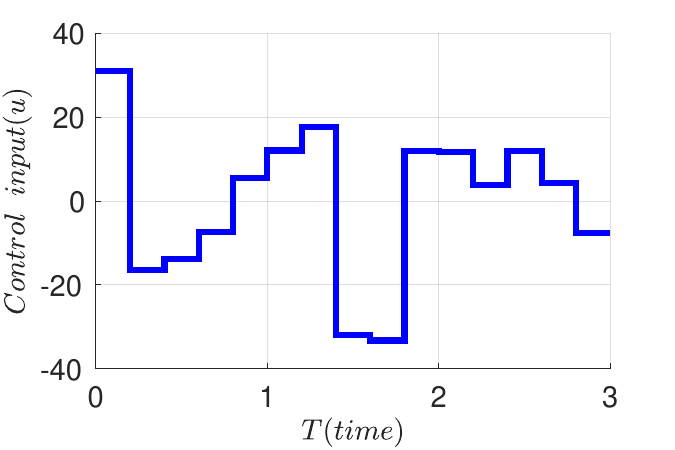}}
    \caption{Optimal control inputs generated through Sequential Quadratic Programming (SQP) optimization using a step size of 0.2 for the cartpole, demonstrating the influence of discretization and the loss of convexity.}
    \label{fig:SQP_cartpole_0.2}
\end{figure}

\begin{figure}[!htbp]
    \centering
   \subfloat[Solution 1 ($J = 1.815065e+03$).]{\includegraphics[width=.49\linewidth]{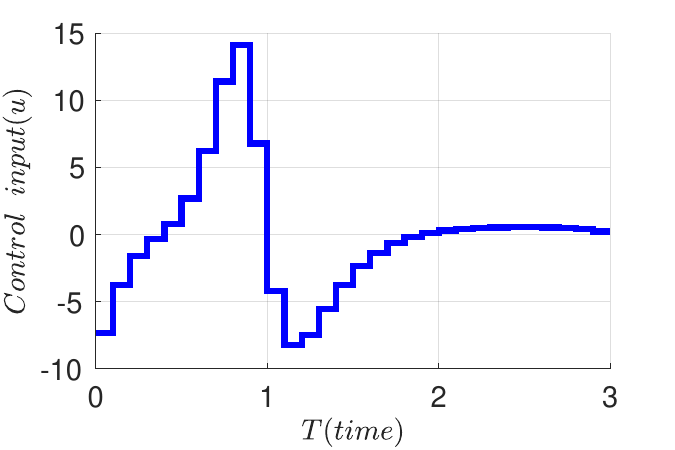}}
   \hfill
    \subfloat[Solution 2 ($J = 7.792933e+03$).] {\includegraphics[width=0.49\linewidth]{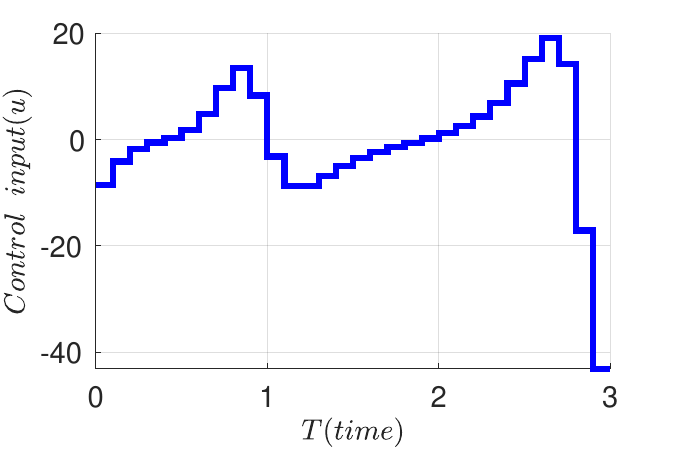}}%
    \hfill
    \subfloat[Solution 3 ($J = 6.108920e+03$).] {\includegraphics[width=0.49\linewidth]{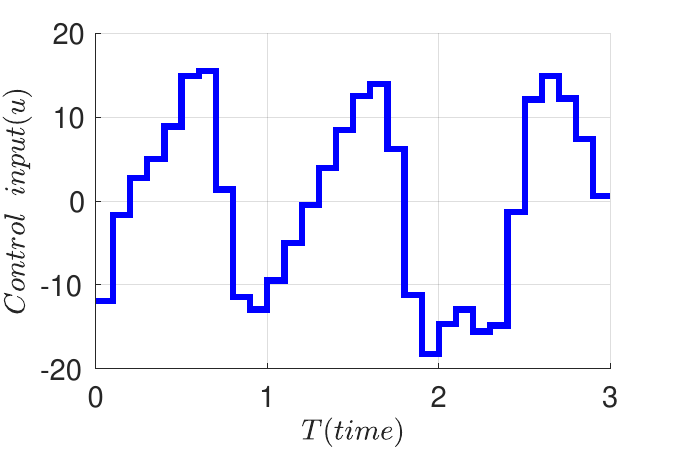}}
    \hfill
    \subfloat[Solution 4 ($J = 1.939614e+04$).] {\includegraphics[width=0.49\linewidth]{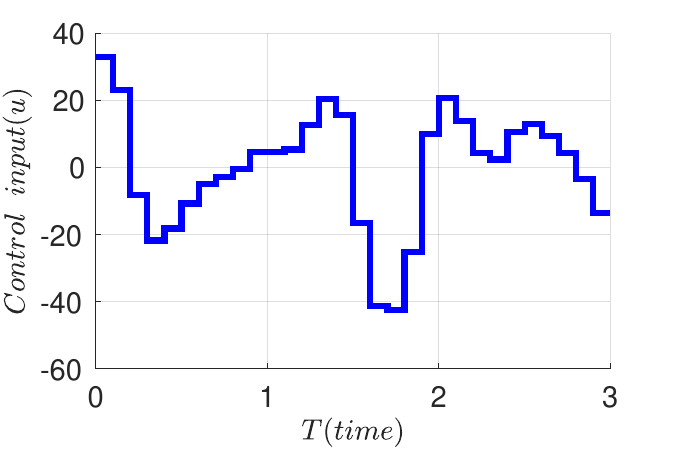}}
    \caption{Optimal control inputs generated through Sequential Quadratic Programming (SQP) optimization using a step size of 0.1 for the cartpole, demonstrating the influence of discretization and the loss of convexity.}
    \label{fig:SQP_cartpole_0.1}
\end{figure}

In this subsection, we set up the discrete-time problem. The problem is solved in MATLAB using \textit{fmincon}. The discrete-time cost becomes 
\begin{align}
    J = \frac{1}{2}\tr{x_{N}}S_{f}x_{N} + \Delta t \sum_{i = 0}^{N-1} \tr{x_{i}}Qx_{i} + \tr{u_{i}}Ru_{i},  
\end{align}
where the subscript $N$ denotes the terminal state. We feed the dynamics to \textit{fmincon} as nonlinear constraints. To do this, we are propagating the dynamics using Euler transcription.\\
We begin by simulating the system for $\Delta t$ = $0.2$. We start with random initial guesses on the control and state variables. We have drawn these random guesses from a normal distribution with zero mean. Figure (\ref{fig:SQP_cartpole_0.2}) shows the result obtained from SQP for the discretized system. We have presented four of the many optimal solutions obtained for this discretization. It is apparent that when using a discretization of $\Delta t = 0.2$, the optimal control problem yields multiple solutions. Additionally, it should be kept in mind that there may exist other spurious solutions, although it is very difficult to find them all. Next, we refined the discretization by selecting a time step of $\Delta t = 0.1$. The results obtained are shown in Fig. (\ref{fig:SQP_cartpole_0.1}). Once more, we observed convergence to multiple spurious solutions showing that the time discretization is still subpar. 

\begin{figure}[!htbp]
    \centering
   \subfloat[Solution 1 ($J = 1.785849e+03$).]{\includegraphics[width=.49\linewidth]{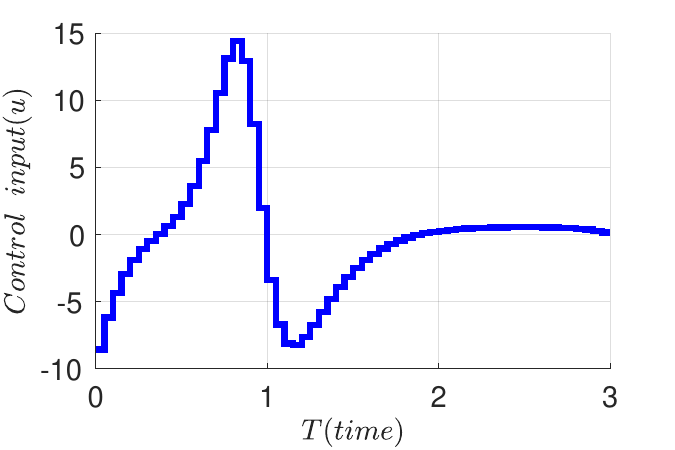}}
   \hfill
    \subfloat[Solution 2 ($J = 5.979149e+03$).] {\includegraphics[width=0.49\linewidth]{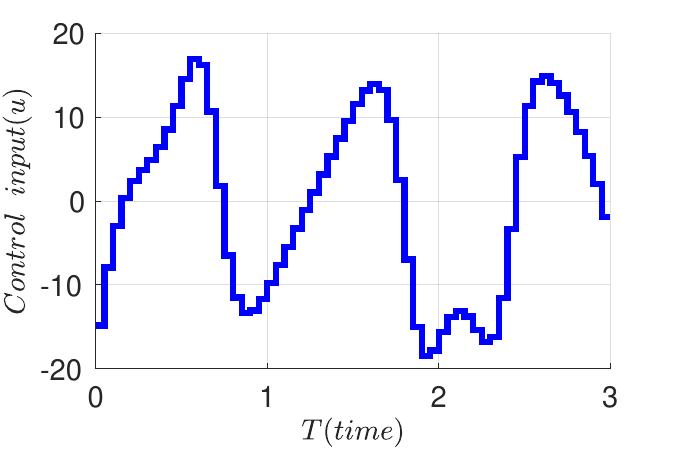}}%
    \caption{Optimal control inputs generated through Sequential Quadratic Programming (SQP) optimization using a step size of 0.05 for the cartpole, demonstrating the influence of discretization and the loss of convexity.}
    \label{fig:SQP_cartpole_0.05}
\end{figure}

\begin{figure}[!htbp]
    \centering
   {\includegraphics[width=0.8\linewidth]{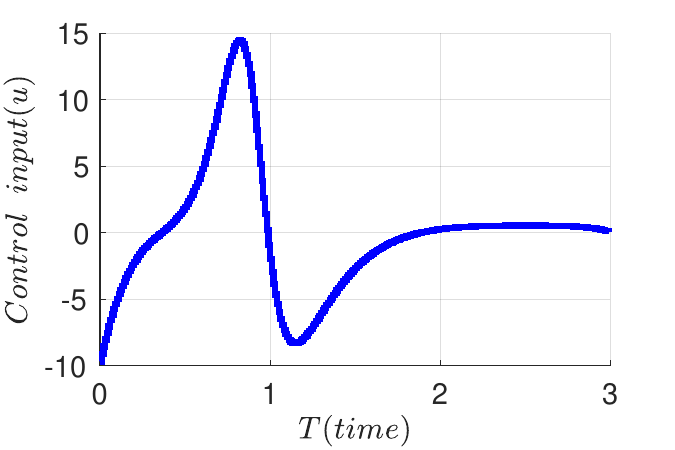}}
   \caption{Unique optimal control trajectory obtained from SQP for Cartpole with discretization step $\Delta t$ = 0.01 ($J = 1.771182e+03$).}
   \label{fig:SQP_cartpole_0.01}
\end{figure}

\begin{figure}[!htbp]
    \centering
    {\includegraphics[width=1\linewidth]{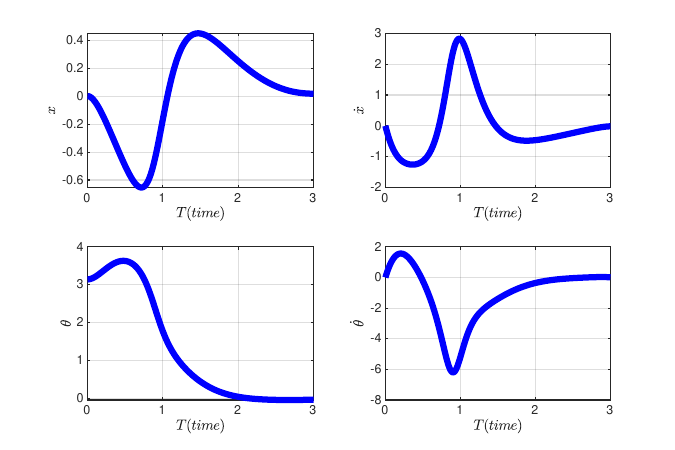}}
   \caption{Evolution of states for the unique solution obtained from SQP for Cartpole with discretization step $\Delta t$ = 0.01.}
   \label{fig:SQP_final_states}
\end{figure}

\begin{figure}[!htbp]
    \centering
   {\includegraphics[width=0.8\linewidth]{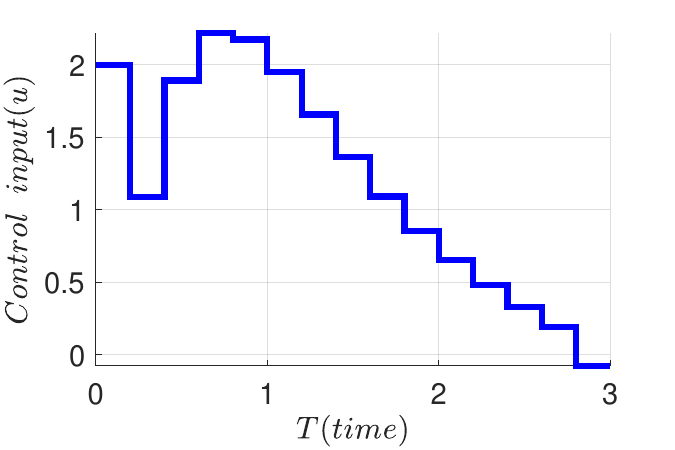}}
   \caption{Optimal control solution obtained from SQP for pendulum with discretization step $\Delta t$ = 0.2.}
   \label{fig:sqp_pen}
\end{figure}

\begin{figure}[!htbp]
    \centering
   {\includegraphics[width=.49\linewidth]{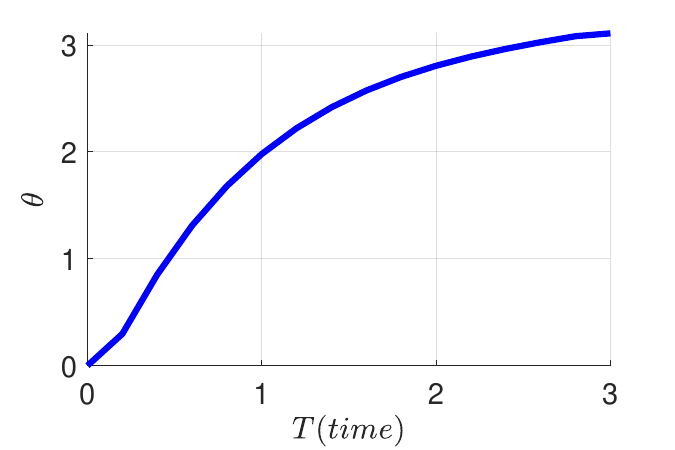}}
   \hfill
     {\includegraphics[width=0.49\linewidth]{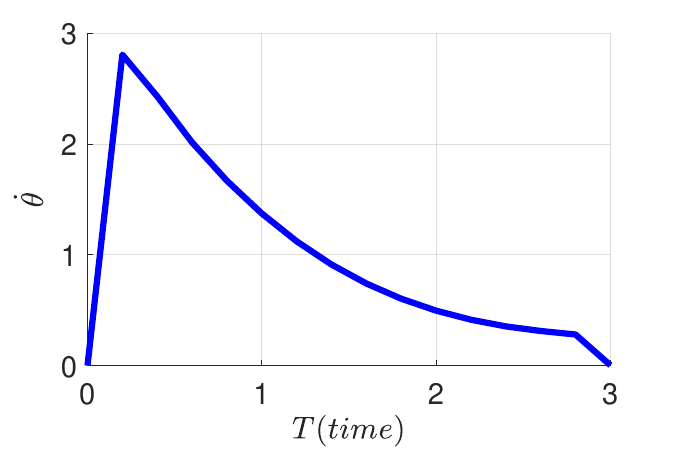}}%
    \caption{Evolution of states for unique control trajectory obtained from SQP for Pen with discretization step $\Delta t$ = 0.2.}
    \label{fig:sqp_pen_states}
\end{figure}

This validates our theory that if the discretization is not fine enough, we will not have a unique solution. In search of the unique solution, we further reduce the step size to smaller values until we converge on a single solution for every simulation. For the next set of simulations, we reduce the time step to $\Delta t$ = 0.05. Figure (\ref{fig:SQP_cartpole_0.05}) shows the results obtained for this discretization. Even for this discretization, the solution does not converge at the same value. Finally, when we reduce the discretization to $\Delta t = 0.01$, we obtain one solution for any random initial guess. This validates our theory of the unique solution if the discretization is fine enough.

Similar to the cartpole, we simulate the discrete-time problem for the pendulum. We started with a time-step of $\Delta t$ = $0.2$, and the problem converged to the same solution for different random initial guesses. Thus, this time discretization is fine enough for the pendulum. The optimal control input obtained is represented in Fig. (\ref{fig:sqp_pen}). The corresponding evolution of states is given in Fig. (\ref{fig:sqp_pen_states}).

\subsection{A Comparison of SQP with iLQR}
\begin{figure}[!htbp]
    \centering
   \subfloat[Solution obtained from iLQR for cartpole with discretization step $\Delta t = 0.2$ ($J = 1.941845e+03$).]{\includegraphics[width=0.49\linewidth]{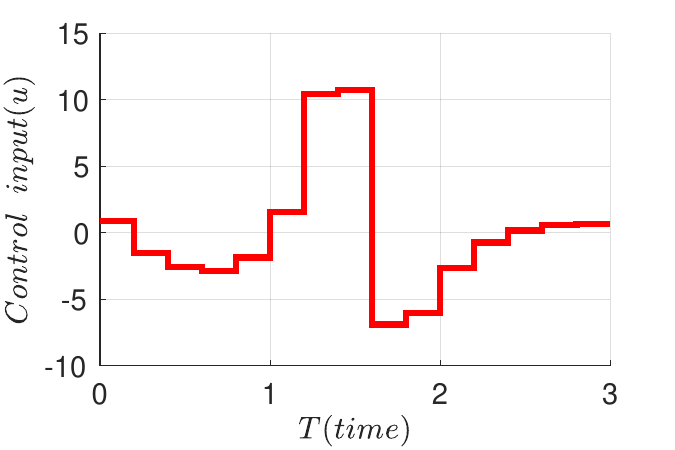}}
   \hfill
   \subfloat[Solution obtained from iLQR for cartpole with discretization step $\Delta t = 0.1$. ($J = 1.815067e+03$).]{\includegraphics[width=0.49\linewidth]{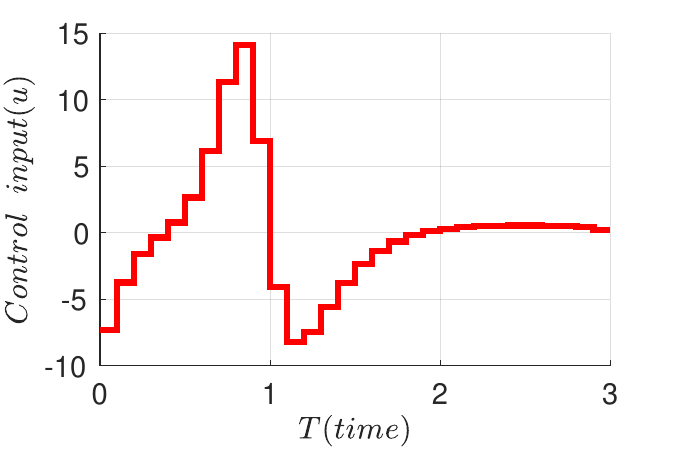}}
   \hfill
   \subfloat[Solution obtained from iLQR for cartpole with discretization step $\Delta t = 0.05$ ($J = 1.785852e+03$).]{\includegraphics[width=0.49\linewidth]{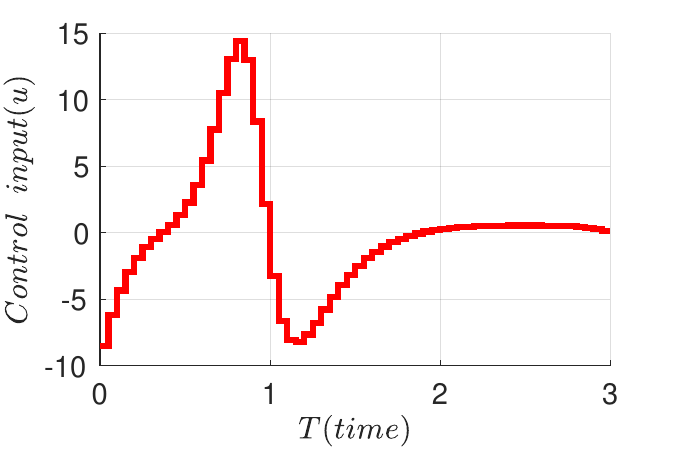}}
   \hfill
   \subfloat[Solution obtained from iLQR for cartpole with discretization step $\Delta t = 0.01$ ($J = 1.771184e+03$).]{\includegraphics[width=0.49\linewidth]{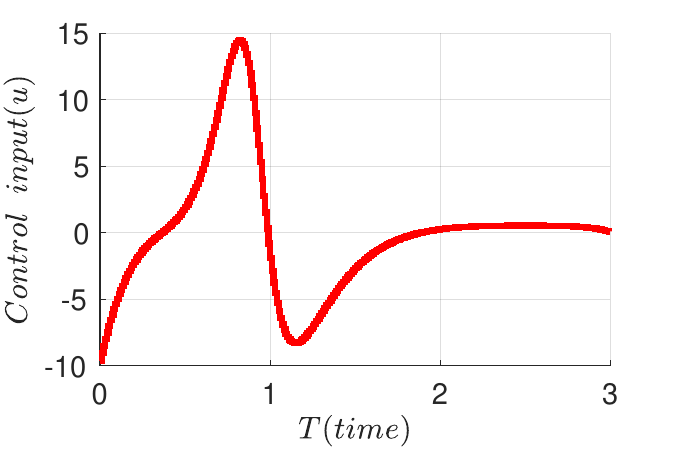}}
   \caption{An illustration of convergence of iLQR for discretization. iLQR converged to the same solution starting from random initial guesses, unlike SQP. }
   \label{fig:ilqr_control_cartpole}
\end{figure}

\begin{figure}[!htbp]
    \centering
   {\includegraphics[width=0.8\linewidth]{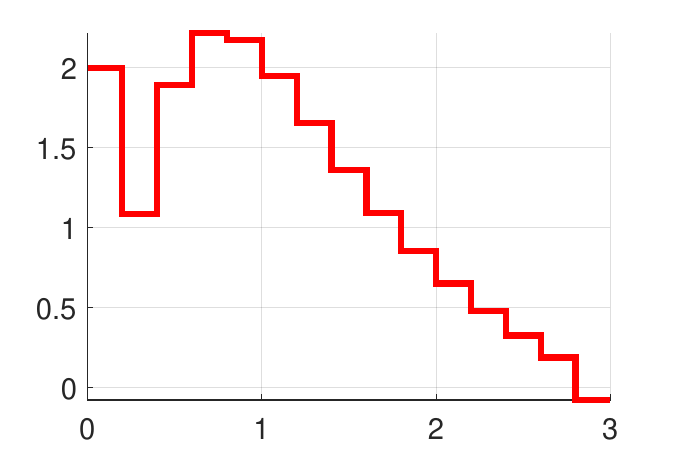}}
   \caption{Optimal control solution obtained from SQP for pendulum with discretization step $\Delta t$ = 0.2.}
   \label{fig:ilqr_pendulum_control}
\end{figure}

To validate the assertions made in section \ref{sec:comp_meth}, we conduct simulations using iLQR for both the pendulum and cartpole systems (see Figure (\ref{fig:ilqr_control_cartpole}) and (\ref{fig:ilqr_pendulum_control})). Maintaining consistency, we retained all parameters identical to those utilized in the SQP simulations to ensure comparability of results. Notably, unlike SQP, iLQR consistently converged to the same solution despite using random initial guesses drawn from a normal distribution with zero mean. We note that this property was also noticed in the original ILQR paper \cite{ilqg} and thought to signify that the ILQR always converges to the global minimum.

\begin{figure}[!htbp]
    \centering
   {\includegraphics[width=1\linewidth]{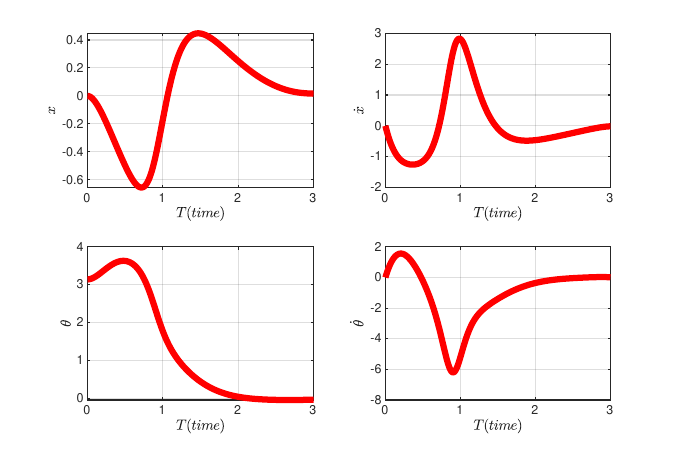}}
   \caption{Evolution of states for unique control trajectory obtained from iLQR for Pen with discretization step $\Delta t$ = 0.01.}
\end{figure}
\begin{figure}[!htbp]
    \centering
    {\includegraphics[width=0.49\linewidth]{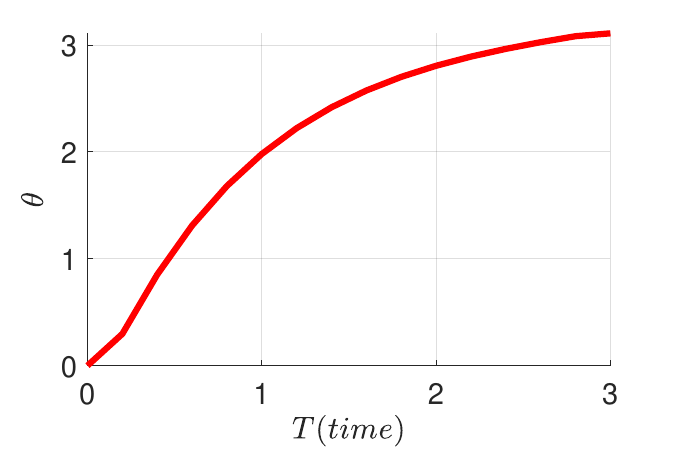}}
   \hfill
   {\includegraphics[width=0.49\linewidth]{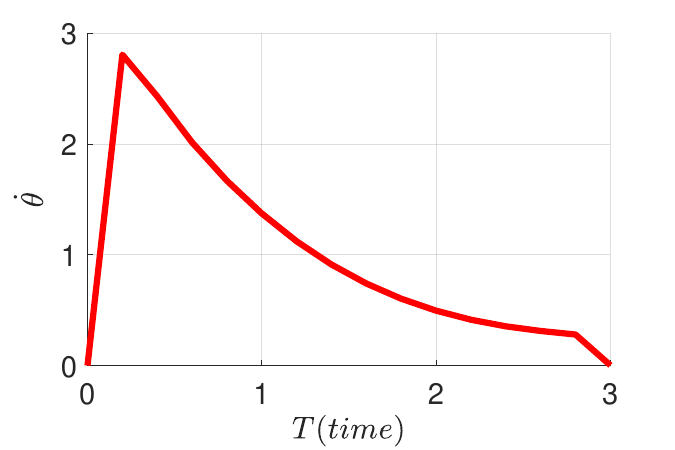}}
   \caption{Evolution of states for unique control trajectory obtained from SQP for pendulum with discretization step $\Delta t$ = 0.2.}
   \label{fig:ilqr_pen_states}
\end{figure}

Thus, we investigate the behaviour of iLQR by initializing it with alternative minima obtained from SQP optimization. Interestingly, iLQR gets stuck at these minima, confirming their validity as solutions that indeed satisfy the necessary conditions of optimality. However, when iLQR is initialized with random initial conditions rather than these spurious minima, it tends to converge towards the ``global minima." This observation is supported by Fig. (\ref{fig:SQP_cartpole_0.2}), where the minimum cost curve aligns with the convergence of iLQR shown in Fig. (\ref{fig:ilqr_control_cartpole}) for the same discretization. \\
For the pendulum, the results from iLQR are represented in Figs. (\ref{fig:ilqr_pendulum_control}) and (\ref{fig:ilqr_pen_states}), which is the same as results obtained from SQP in Figs. (\ref{fig:sqp_pen}) and (\ref{fig:sqp_pen_states}). We have the same solution from SQP and iLQR as the initial discretization of $\Delta t$ = 0.2 is good enough to maintain the convexity of the Hamiltonian and give a unique solution.\\
During our simulation, we observed that beyond demonstrating superior convergence performance, iLQR exhibits lower time complexity compared to the NLP approach and, hence, is computationally far more efficient.
Thus, it is imperative to emphasize the significance of selecting an appropriate time discretization when solving an optimal control problem, as an inadequately chosen time step may compromise the structural integrity of the problem, leading to spurious minima that are not present in the continuous problem. 

\section{{Conclusion}}\label{section:conclusion}
In this paper, we investigate continuous-time optimal control problems distinguished by the unique attribute of possessing a convex Hamiltonian. We provide a theoretical analysis demonstrating the uniqueness of solutions for such problems given the strict convexity of the Hamiltonian in the control variable. However, upon discretizing the dynamics, we observe the emergence of multiple solutions as the convexity is compromised. To empirically validate the theory, we set up a quadratic cost discrete-time optimal control problem for the cartpole and pendulum swing-up tasks and solved it using SQP and iLQR. SQP converged to various different minima for coarse discretization and as the discretization is made finer, the SQP algorithm converges to the unique solution for any initial guess. Notably, despite the existence of multiple solutions in the discrete-time scenarios, iLQR consistently converges to the same solution across various initializations, owing to the intrinsic forward-backwards structure of the algorithm. Additionally, we note that iLQR consistently converges to the minimum cost solution, even in cases of multiple solutions resulting from coarse discretization. In future work, we need to investigate the theoretical reason for the convergence of ILQR to the ``global minimum" regardless of the time discretization and also be able to extend the analysis to cases with state and control constraints.

\printbibliography

\end{document}